\documentclass{article}

\usepackage{graphicx}
\usepackage{amsmath}
\usepackage{color}
\usepackage{amssymb}
\usepackage{epsfig}
\usepackage{psfrag}

\newcommand{\ds}{\displaystyle}

\newcommand{\pr}{\partial}
\newcommand{\R}{\mathbb R}
\newcommand{\e}{{\rm e}}
\newcommand{\atan}{\tan^{-1}}
\newcommand{\nfrac}[2]{\mbox{\large $\textstyle \frac{#1}{#2}$}}
\newcommand{\mat}{{\sc Matlab}}

\usepackage{enumerate}

\begin{document}

\title
{Exact solutions for logistic reaction-diffusion in biology}

\author{P. Broadbridge$^1$, B.H. Bradshaw-Hajek$^2$}


\date{\small 1. Dept. of Mathematics and Statistics, La Trobe University, Victoria, Australia. \\
2. Phenomics and Bioinformatics Research Centre, School of Information Technology and Mathematical Sciences, University of South Australia.  Bronwyn.Hajek@unisa.edu.au}

\maketitle

%


\begin{abstract}
\noindent Reaction-diffusion equations with a nonlinear source have been widely used to model various systems, with particular application to biology. Here, we provide a solution technique for these types of equations in $N$-dimensions. The nonclassical symmetry method leads to a single relationship between the nonlinear diffusion coefficient and the nonlinear reaction term; the subsequent solutions for the Kirchhoff variable are exponential in time (either growth or decay) and satisfy the linear Helmholtz equation in space. Example solutions are given in two dimensions for particular parameter sets for both quadratic and cubic reaction terms.
\end{abstract}


\section{Introduction}

Logistic reaction-diffusion equations with a nonlinear source are widely used to model many different systems, particularly in biology. One of the earliest appearances of such a model was in the seminal paper by Fisher in 1937 \cite{F37}, where he introduced the equation
\begin{equation}\label{Fisher}
\ds\frac{\pr \theta}{\pr t}=D\ds\frac{\pr^2\theta}{\pr x^2}+s\theta(1-\theta).
\end{equation}
Fisher's equation originally modelled the frequency in a diploid population, of a new advantageous recessive gene, labelled `$a$'. In fact, for a sexually reproducing species, Fisher's equation follows only if there are three different phenotypes $AA$, $Aa$ and $aa$ whose relative fitness coefficients are in linear progression; otherwise Fisher's assumptions lead to a cubic rather than a quadratic source term \cite{S73,BBFA02}. Fisher's equation remains the model of choice for many biological problems such as those in population dynamics (where $\theta(x,t)$ is the population density divided by the carrying capacity of the environment \cite{Schaefer}), and for biological cellular tissue growth (where $\theta(x,t)$ is the cell population density divided by the steady-state tissue density, see for example \cite{L}).

Since cell mobility depends on cell density, such models naturally generalise to nonlinear reaction-diffusion with logistic reaction and nonlinear diffusion. Extending the model to three dimensions and including a broader range of source terms, the general form of a logistic reaction-diffusion equation for a density $\theta({\bf r},t)$ is
\begin{equation}\label{ge}
\theta_t=\nabla\cdot[D(\theta)\nabla \theta]+R(\theta)
\end{equation}
where $D(\theta)$ is the density dependent diffusion coefficient, $R(\theta)$ is a logistic type source term which may be quadratic or cubic, and $\nabla$ is the usual gradient operator in three dimensions. Some properties of these equations are shared more generally by equations of Kolmogorov-Piscounov-Petrov type \cite{AW75} but here, we consider $R(\theta)$ to be either quadratic (as in the traditional Fisher-type model) or cubic (as in the Huxley and Fitzhugh-Nagumo models). The Fitzhugh-Nagumo equation has been used to model a nerve axon potential, the intermediate unstable steady state $\theta=\theta_1$ being the threshold electrical potential that separates the basins of attraction for the stable steady activated state $\theta=1$ and the quiescent state $\theta=0$.

Despite the wide use of equation \eqref{ge}, few exact analytic solutions are known even in the one-dimensional case.  We have not been able to find previously published exact solutions for logistic reaction-diffusion when $D(\theta)$ is non-constant. When $\theta({\rm {\bf r}},t)=\theta(x,t)$ (i.e. the one-dimensional case), $D$ is constant and $R(\theta)$ is quadratic, equation \eqref{ge} is known as the Fisher, Fisher-Kolmogorov or KPP equation due to the classic papers by these authors \cite{F37,KPP37}. An exact travelling wave solution was first presented by Ablowitz and Zeppetella in 1979 \cite{AZ79}. 
In two dimensions, explicit and approximate travelling wave solutions have been presented by Brazhnik and Tyson \cite{BT99}.

When $\theta({\rm r},t)=\theta(x,t)$, $D$ is constant, and $R(\theta)$ is cubic or of higher order Huxley type,   Kametaka \cite{Kametaka} found a travelling wave solution, while other travelling wave solutions have been found by McKean \cite{M70} and Rinzel \cite{R75}. Periodic solutions were found by Carpenter \cite{C77} and Hastings \cite{H74}. Arrigo, et al. \cite{AHB94}, and Clarkson and Mansfield \cite{CM94} found some exact solutions using the method of nonclassical symmetry analysis and these solutions have been applied to the problem of a new advancing recessive gene \cite{BBFA02,BB04}. The solutions found using the nonclassical symmetry approach can also be found using the Painlev\'{e} approach \cite{Conte88,CG92}. Other solutions have been presented by Kawahara and Tanaka \cite{KT83}, Kudryashov \cite{K92}, Chen and Gu \cite{CG99} and Nikitin and Barannyk \cite{NB04}. The existence of solutions in the cubic case was investigated by Nagylaki \cite{N75} and Conley \cite{C75}.

In this paper, we use the nonclassical symmetry method to present exact analytical solutions to equation \eqref{ge} when $R(\theta)$ is quadratic with two real roots, cubic with one doubly repeated root, and cubic with three distinct roots. In each case, the nonlinear diffusion coefficient takes on a particular form.

In Section 2, we describe the solution technique; the solutions are presented in Sections 3, 4 and 5. The application of the solutions to population genetics is discussed in Section 6 and some final remarks are presented in Section 7.

\section{Nonclassical reduction to the Helmholtz equation}

A full Lie point symmetry classification of equation \eqref{ge} was made by Dorodnitsyn et al. \cite{Dorod}. In classical Lie point symmetry analysis, one seeks transformations that leave the governing equation invariant. In some cases, these transformations may be used to simplify the governing equation, leading to a possible analytic solution. Classically invariant solutions include travelling waves and scale-invariant solutions that may exhibit extinction, single-peak or multi-peak blow-up in finite or infinite time, with unbounded or compact spatial support on $\Re^n$, only when the free functions for diffusivity and/or reaction are power-laws, exponentials or of $\theta\log\theta$ form.

Following Ovsiannikov's general formulation of partial symmetries \cite{Ovsiannikov}, the idea of nonclassical symmetries was pioneered by Bluman and Cole \cite{BC69}. This requires a transformation to leave the system consisting of the governing equation {\it and} the invariant surface condition, invariant. This extra requirement can sometimes give rise to transformations that cannot be found by the Lie point method. Nonclassical symmetry methods, also known as Q-conditional symmetries, have been applied to equations belonging to the class \eqref{ge}, and a number of forms of $D(\theta)$ and $R(\theta)$ have been found to admit strictly nonclassical symmetries \cite{AHB94,CM94,AH95,GB96}. The complete nonclassical symmetry classification of equation \eqref{ge} in two dimensions was given by Goard and Broadbridge \cite{GB96}. Some of the same nonclassical symmetries readily extend to $N$-dimensions \cite{BBT15}.

By writing equation \eqref{ge} in terms of the Kirchhoff variable (see for example \cite{P69})
\begin{equation}\label{Kirch}
u=u_0+\ds\int_{\theta_0}^{\theta}D(\theta')~d\theta',
\end{equation}
so that a boundary condition $\theta=\theta_0$ corresponds to $u=0$, we obtain
\begin{equation}\label{linD}
F(u)\ds\frac{\pr u}{\pr t}=\nabla^2u+Q(u)
\end{equation}
where $F(u)=1/D(\theta)$ and $Q(u)=R(\theta)$. This equation admits the nonclassical reduction operator \cite{GB96}
\begin{equation}\label{symgen}
\Gamma=\ds\frac{\pr}{\pr t}+Au\frac{\pr}{\pr u}
\end{equation}
whenever $F$ and $Q$ are related by
\begin{equation}\label{FQrel}
Q(u)=AuF(u)+\kappa u
\end{equation}
for $A,~\kappa\in\R$ constant.  Given $D(\theta)$ this gives $R(\theta)$ by direct integration:
\[
R(\theta)=\ds\frac{A}{D(\theta)}\int  D~d\theta+\kappa\int D d\theta.
\]
However, given $R(\theta)$, $D(\theta)$ is obtained by solving a differential equation:
$$D(\theta)=\frac{du}{d\theta}=\frac{Au}{R-\kappa u},$$
equivalently
\begin{equation}
D({\theta})=\frac{ A\int_{\theta_0}^\theta Dd\theta}{R(\theta)-\kappa\int_{\theta_0}^{\theta} D~d\theta}.
\label{use}
\end{equation}

Equation \eqref{symgen} is a genuine nonclassical symmetry because it leaves equation \eqref{ge} invariant only if one also makes use of the invariant surface condition, $u_t=Au$. However this conditional invariance allows a consistent reduction of the original PDE, to a differential equation among the invariants of the symmetry. Making use of this nonclassical reduction, equation \eqref{linD} can be transformed to the linear Helmholtz equation
\begin{equation}\label{Helm}
\nabla^2\Phi+\kappa\Phi=0 \quad\quad {\rm with} \quad\quad u=\e^{At}\Phi({\bf x}).
\end{equation}
In this manner, an arbitrary solution of the linear Helmholtz equation in $N=1,~2$ or 3 dimensions may be used to construct a solution of the nonlinear reaction-diffusion equation. For example, one may use any of the solutions that have previously been constructed for the amplitude of a scattered acoustic wave \cite{Colton}. However, such solutions that represent scattering by a finite body, must approach the isotropic solution at large distances. Therefore the radial solutions $u(r,t)$ are considered to be canonical, and indicative of features of fields scattered from aspherical boundaries. In the current context, there are three possible types of solution to be considered, namely those with $\kappa=0$, $\kappa<0$ and $\kappa>0$.

  With $\kappa=0$, $\Phi (r)$ is a radial solution of the Laplace equation, which can only be a linear combination of the trivial constant solution and the unit point source solution. Any such non-constant solution must be singular at the origin, where there is a steady flux from or to a point source or sink. The differential equation (\ref{use}) is now linear, with general solution
\begin{equation}
u=c_1\exp\int \frac{A}{R(\theta)}d\theta.
\label{solK0}
\end{equation}

  With $\kappa =-K^2<0$, in two dimensions $\Phi$ must be a linear combination of modified Bessel functions $I_0(Kr)$ and $K_0(Kr)$ which must either be infinite at the origin or be unbounded at large $r$.  A similar situation pertains in three dimensions when the modified Bessel functions are replaced by spherical modified Bessel functions.
  
   Finally, in the case $\kappa=K^2>0$, there are positive bounded solutions that are Bessel functions $\Phi=J_0(Kr)$ in two dimensions and spherical Bessel functions $\Phi=j_0(Kr)$ in three dimensions. They satisfy boundary conditions
 \begin{eqnarray}
 u_r(0,t)=0,\\
 u(r_1,t)=0,
 \end{eqnarray}
 where $r_1=\lambda_1/K$, $\lambda_1$ being the first zero of the Bessel function. This outer boundary condition may represent extreme total culling of a species at a boundary, for example extreme harvesting of a prey species at the boundary of a protection zone, selective culling of some genotype, or removal of outer tumour cells by radiation or chemo-therapy or diathermy. Alternatively, since $p=-u_r(r,t)/u(r,t)=\Phi_r(r)/\Phi(r)$ varies from 0 at $r=0$ to $\infty$ at $r=r_1$, the outer boundary may be relocated to some location $r_2\in (0,r_1)$ where it satisfies a Robin condition $-u_r=pu$ for some pre-chosen parameter $p$. Since $-u_r$ is simply the radial flux density of the population, this may approximately represent individuals responding to an external chemoattractant with a fixed probability proportional to $p$, except that when the diffusivity is not constant, $u$ is not exactly proportional to the population $\theta$.

For convenience, from here on we set $\theta_0=0$. 

A solution $D(\theta)$ of (\ref{use}) must be a fixed point of the map
\begin{equation}
D_{n+1}(\theta)=\ds\frac{A\int_0^{\theta} D_{n}~d\theta}{R(\theta)-\kappa\int_0^{\theta} D_n~d\theta}.
\label{Dmap}
\end{equation}
For some values of the system parameters, this is a contraction map that converges to a unique solution \cite{BBT15}.
The modelling constraint $D(\theta)>0$ may restrict the values of the temporal exponential  growth parameter $A$ that can occur. These details will depend on the form of the reaction function $R(\theta)$.
We now present some exact analytic solutions in the case where $R(\theta)$ is quadratic or cubic.

\section{Fisher-type logistic equations, $R(\theta)=s\theta(1-\theta)$}

The standard generalisation of Fisher's equation to $N$-dimensions is
$$\frac{\partial \theta}{\partial t}=D\nabla^2u+s\theta(1-\theta).$$
Since cell mobility depends on cell density, Fisher's equation naturally generalises to nonlinear reaction-diffusion \eqref{ge} with nonlinear diffusion as well as logistic reaction. Even in one spatial dimension, there are very few known exact solutions, apart from the one-dimensional travelling wave solution with a special non-minimal group velocity \cite{AZ79}. The logistic source changes sign at the carrying capacity, quite different from the positive definite source term of combustion. Despite that major difference, the construction that was previously applied to combustion modelling \cite{BBT15} still can be applied to population dynamics after some restrictions on the system parameters.

	First consider the case $\kappa=0$. The diffusivity is given explicitly by
$$D=-\frac{A}{s}\theta^{-2}\left(\theta^{-1}-1\right)^{A/s}.$$
Since $D(\theta)$ must be positive for $\theta>0$, it must be true that $A<0$. Then $D(\theta)\to\infty$ as $\theta\to 1^-$. Since in biological applications the diffusivity must be bounded, the case $\kappa=0$ is inadmissible.

	Secondly, consider the case $\kappa=-K^2<0$. Since $R(1)=0$, equation (\ref{use}) then implies $A=K^2 D(1)>0$. This can represent only a growing population (by equation \eqref{Helm}). If we presume that $D$ has an upper bound, then exponential growth of $u$ implies unbounded growth of the population density $\theta$. However, this is problematic for population modelling since the source term is negative when $\theta>1$. The radial solution is a linear combination of modified Bessel functions that have a point source or a singularity at the origin. Alternatively, $D(\theta)$ may diverge at some finite value of $\theta$ so that $\theta$ remains bounded as $u$ diverges. However, an unbounded diffusivity is untenable in population modelling.

Finally, we consider the applicable case $\kappa=K^2>0$. It follows from (\ref{use}) that $D(1)=-A/K^2$. Therefore the only valid solutions remaining are those with $A<0$, so that $u$ approaches zero exponentially in time (by equation \eqref{Helm}).
Since $R(\theta)$ is now analytic at $\theta=0$ with leading order $\mathcal{O}(\theta^1)$, compared to the case of Arrhenius combustion \cite{BBT15}, $R(\theta)$ is no longer negligible compared to $u(\theta)$ at small $\theta$. From (\ref{use}) as $\theta\to 0$,
$$D(0)=\lim_{\theta\to 0}\frac{AD(0)\theta}{s\theta-\kappa D(0)\theta}=\frac{AD(0)}{s-\kappa D(0)},$$
implying \begin{equation}
A=s-\kappa D(0).
\label{Aeq}
\end{equation}
For example if $r_1$ is the radius of a circular aquatic reserve, outside of which the population of a mobile species is practically zero, one may assume $\theta=0$ at $r=r_1$ from which it follows that $K=\lambda_1/r_1$ ($\lambda_1$ being the first zero of Bessel function $J_0$), and (\ref{Aeq}) gives a condition that guarantees the non-existence of the undesirable solution with $A<0$, that approaches extinction. This condition may be expressed
\begin{equation}
r_1>\lambda_1\sqrt{D(0)/s}.
\end{equation}
For example, for a mobile species in which an individual's range expands to $100$ km$^2$ per year and the time scale for uninhibited exponential growth is ${1}/{s}=5$ years, the conservative safe diameter of a marine park would be $2r_1\approx108$ km.

\begin{figure}
\includegraphics[trim = 35mm 100mm 35mm 100mm, clip, width=100mm]{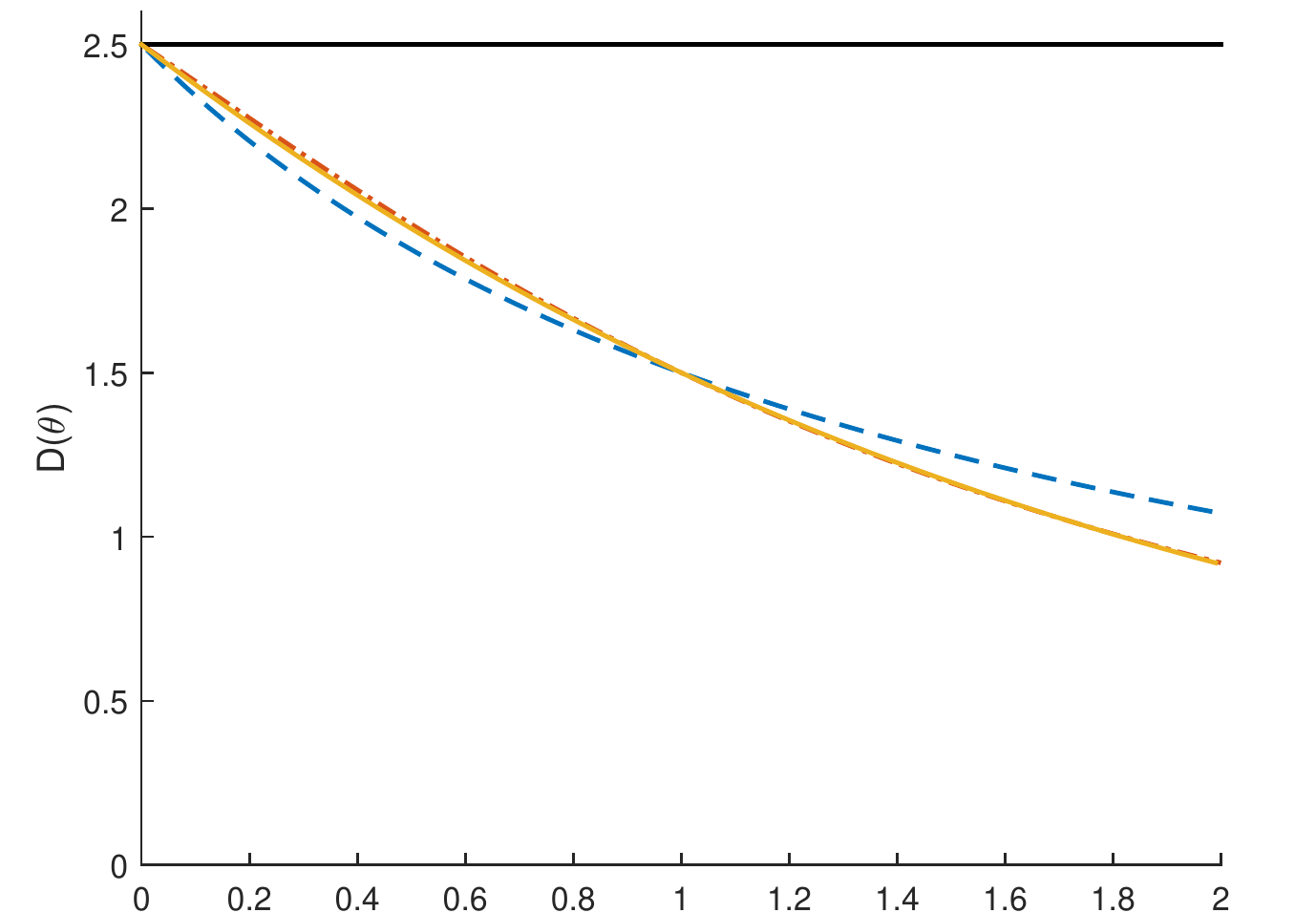}
\caption{$D(\theta)$ for the Fisher source term constructed by \mat~  routine ode45 (solid), as well as approximations $D_0$ (solid), $D_1(\theta)$ (dashed), and $D_2(\theta)$ (dash-dot) ($K=1$, $s=1$ and $A=-1.5$). }
\label{fig:2} 
\end{figure}

By choosing as the first estimate for $D(\theta)$, the constant value $D_0=D(0)=(s+|A|)/K^2$, it follows that all subsequent iterates $D_j(\theta)$ of the map (\ref{Dmap}), must have the correct values of $D(\theta)$ at $\theta=0$ and $\theta=1$. The first few iterates are:
\begin{eqnarray}
D_0&=&D(0),\\
D_1&=&\frac{|A|D_0/s}{\theta+K^2D_0-1},\\
D_2&=&\frac{-A^2D_0\log([s\theta+K^2D_0-s]/[K^2D_0-s])}{s^2\theta(1-\theta)+K^2AD_0\log([s\theta+K^2D_0-s]/[K^2D_0-s])}.
\end{eqnarray}
These are shown in Figure \ref{fig:2} for the example with $s=1$, $K=1$ and $A=-1.5$. 
The iterated approximation $D_2$ closely agrees with the numerical approximation obtained by solving (\ref{use}) with the third/fourth-order Runge-Kutta routine ode45 of \mat. From $\theta=0$ to $\theta=2$ which is well over double the carrying capacity ($\theta=1$), the diffusivity is positive-valued and decreasing. 

\begin{figure}
\includegraphics[trim = 35mm 100mm 35mm 100mm, clip, width=100mm]{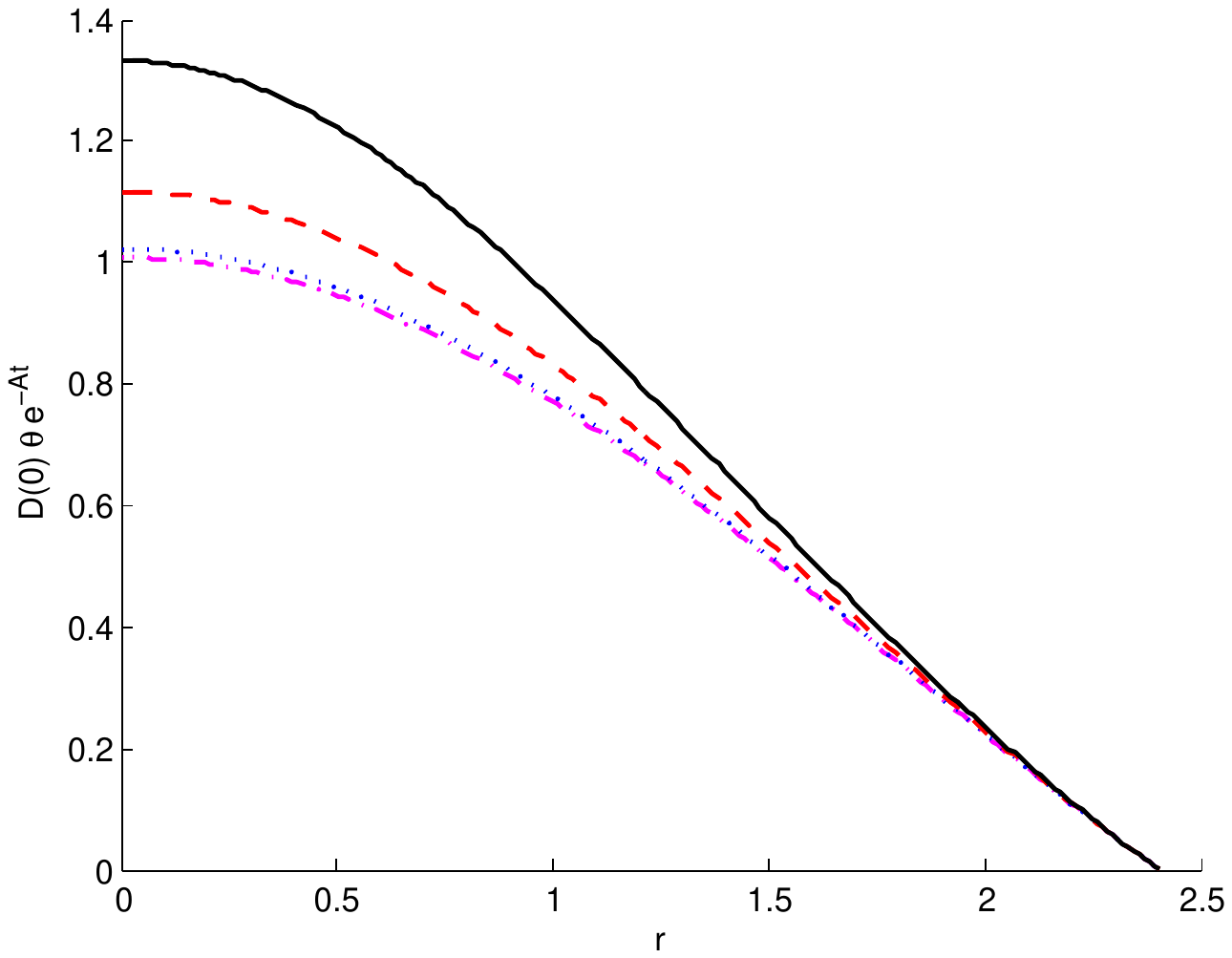}
\caption{Solution for the decreasing population density versus radial coordinate, at times $|A|t=-1.5, 0.0, 1.5, 2.5$, according to the Fisher-type equation with $K=1$ and $s=1$.}
\label{fig:1} 
\end{figure}

The exact solution for the doomed population is shown in Figure \ref{fig:1}.
The shape of the population curve approaches that of the Kirchhoff variable but from above in this case, since $D(\theta)$ is now a decreasing function of population density, a reasonable model for some populations and for cells whose mobility is decreased by close packing.

\section{Huxley equation, with $R(\theta)=s\theta^2(1-\theta)$}

In this case, the equation in $N$ dimensions with nonlinear diffusivity may be written
\[
\theta_t=\nabla\cdot[D(\theta)\nabla \theta]+s\theta^2(1-\theta).
\]
This better models the frequency of a new advantageous recessive gene in the important case considered by Fisher, that of a Mendelian diploid sexually reproducing population. It also occurs in the Huxley approach to modelling the electrical potential in a nervous system.

We first consider the case $\kappa=0$. Equation (\ref{use}) has general solution
\[
D(\theta)=c_1\theta^{-1}(1-\theta)^{-2}(1-\theta^{-1})^{-(A/s)+1}\exp\left(-\frac{A}{s\theta}\right)
\]
This expression is divergent as $\theta\to0$ for $A\ge0$, and also divergent as $\theta\to1^-$ when $A<0$. As such, the case when $\kappa=0$ is inadmissible for modelling a population.

Secondly, we consider the case when $\kappa=-K^2<0$. Since $R(1)=0$, equation \eqref{use} implies that $A=K^2D(1)>0$. As described in the previous section, the case when $A>0$ is not of interest in population modelling.

Finally, consider the case when $\kappa=K^2>0$. Since $R(1)=0$, we find that $A=-K^2D(1)$, so that any valid solutions have $A<0$, and by equation \eqref{Helm}, $u$ approaches zero exponentially in time. When $\theta$ is small, $R(\theta)\sim s\theta^2$, so that by using a Taylor series approximation for $D(\theta)$, in the limit as $\theta\to0$,
\[
D(0)=\lim_{\theta\to0}\ds\frac{AD(0)\theta}{s\theta^2-K^2D(0)\theta}=-\ds\frac{A}{K^2},
\]
implying $A=-K^2D(0)$.

For a problem in 2-dimensions, the solution to equation \eqref{Helm} is the Bessel function $J_0(Kr)$. Again, if we assume, for example, that $r_1$ is the radius of an aquatic reserve and that outside the reserve the population is practically zero, $\theta=0$ at $r=r_1$, then $r_1=\lambda_1/K$, where $\lambda_1$ is the first zero of the Bessel function $J_0$. Within the reserve, $u>0$ and so $\theta>0$.

It would be desirable if we could calculate a minimum radius for the aquatic reserve such that the population inside could be sustained and not become extinct. However, this can only happen if $A\ge0$, and since $A=-K^2D(0)$ and $K\ne0$, this is not possible. As a result, the extinguishing solution will always exist for a population that can be appropriately modelled using the Huxley source term.

\begin{figure}
\includegraphics[trim = 35mm 100mm 35mm 100mm, clip, width=100mm]{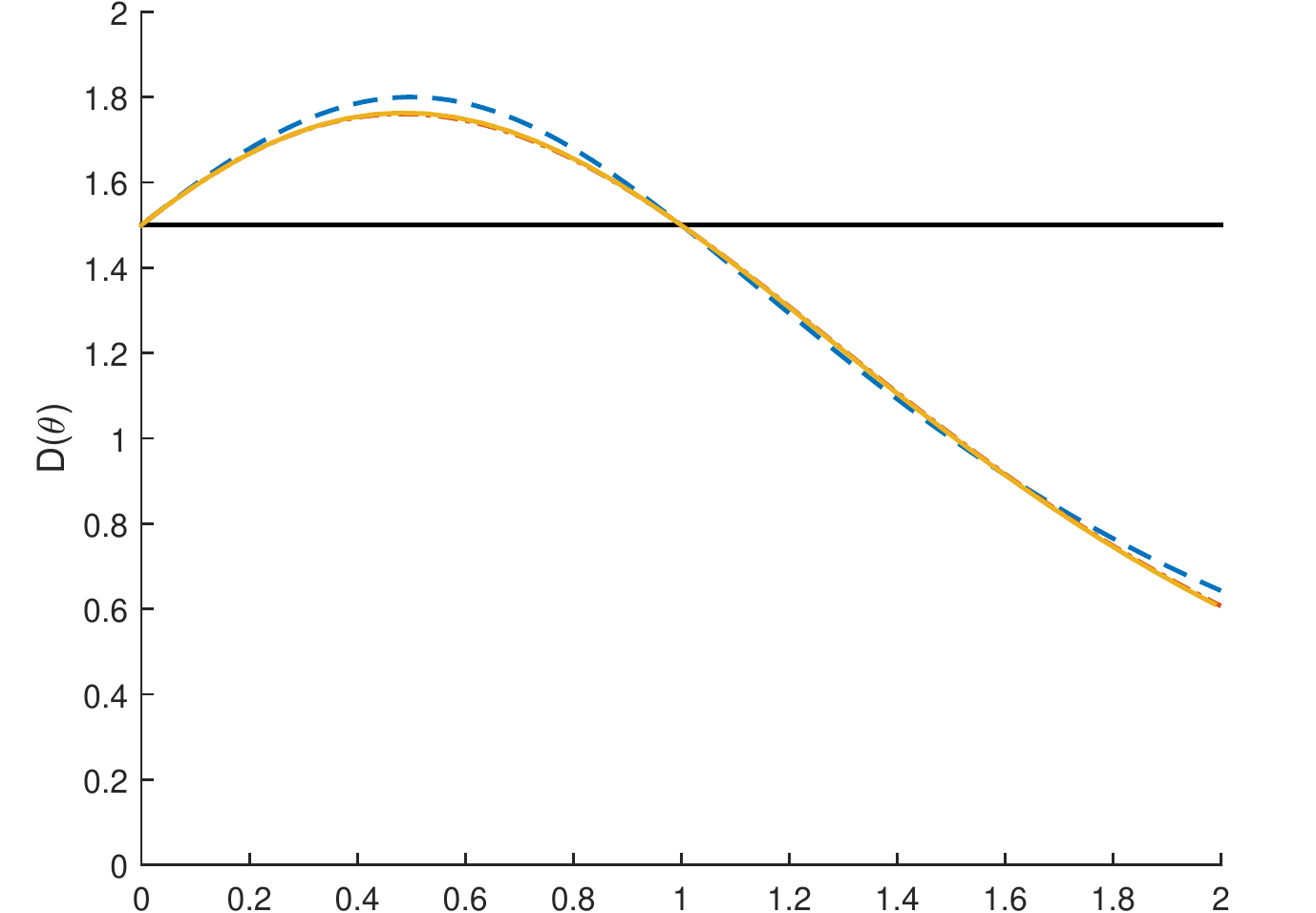}
\caption{$D(\theta)$ for the Huxley source term, constructed by \mat routine ode45 (solid), as well as approximations $D_0$ (solid), $D_1(\theta)$ (dashed), and $D_2(\theta)$ (dash-dot) ($K=1$, $s=1$ and $A=-1.5$). The second iterate, $D_2(\theta)$, is almost indistinguishable from the numerical solution.}
\label{fig:4} 
\end{figure}

The nonlinear diffusivity $D(\theta)$ can be calculated in the same way as for the Fisher case. By choosing $D_0=D(0)=-A/K^2$ for the first estimate for $D(\theta)$ we find that
\begin{eqnarray}
D_0&=&D(0),\\
D_1&=&-\ds\frac{A^2}{K^2}\ds\frac{1}{s\theta(1-\theta)+A}~,\\
D_2&=&\ds\frac{A^3}{K^2}\ds\frac{\atan\left(\frac{1}{\beta}\left(\theta-\nfrac12\right)\right)+\atan\left(\frac{1}{2\beta}\right)} {\beta s^2\theta^2(1-\theta)-A^2\left[\atan\left(\frac{1}{\beta}\left(\theta-\nfrac12\right)\right)+\atan\left(\frac{1}{2\beta}\right)\right]}~,
\end{eqnarray}
where $\beta=\sqrt{\left|\frac{A}{s}+\nfrac14\right|}$. These iterates for the nonlinear diffusion are shown in Figure \ref{fig:4}, together with the solution to equation \eqref{use} found using \mat's ode45 routine.

\begin{figure}
\includegraphics[trim = 35mm 100mm 35mm 100mm, clip, width=100mm]{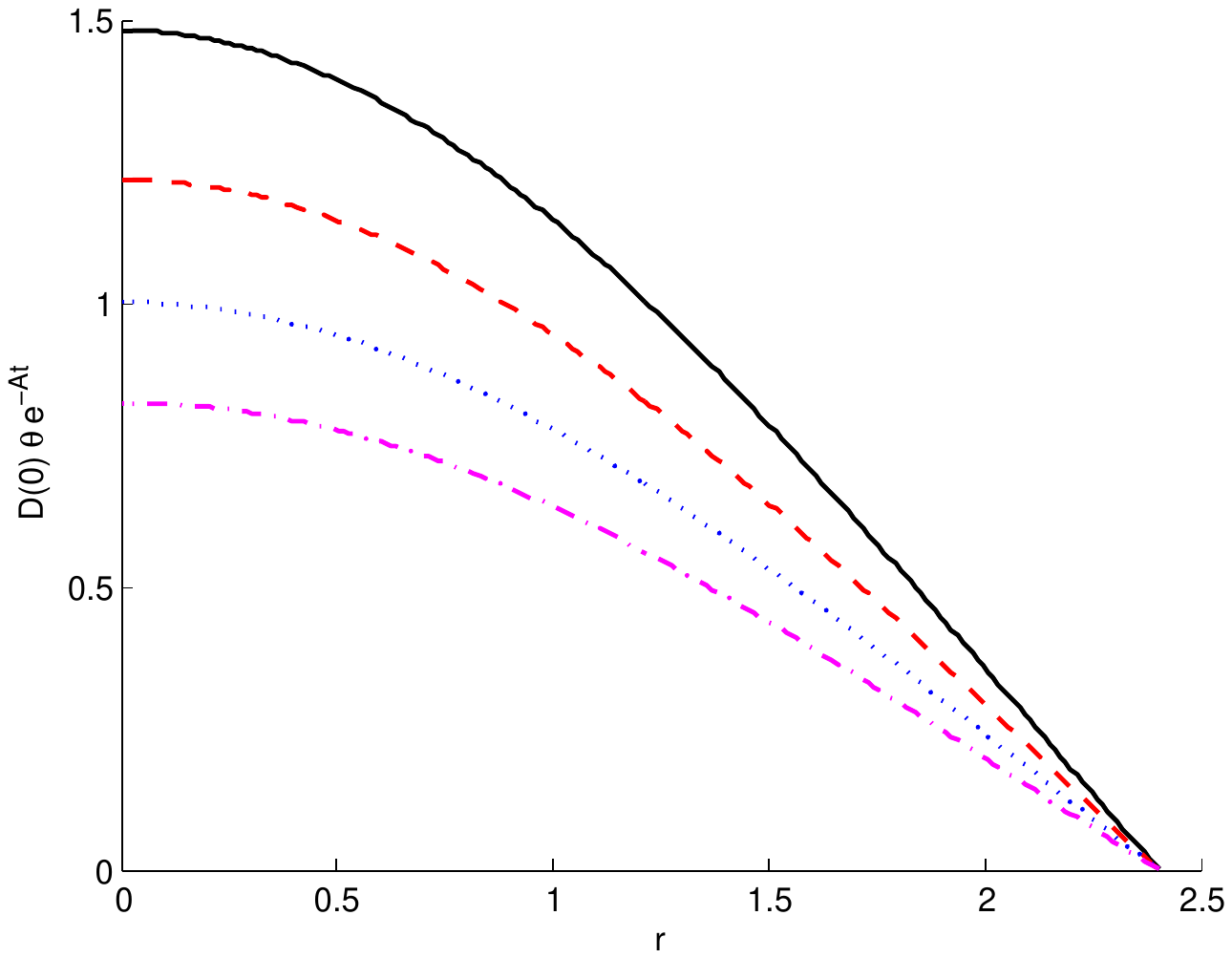}
\caption{Solution for the decreasing population density versus radial coordinate, for the Huxley source term, at times $|A|t=0.0, 0.1, 0.2, 0.3$, according to the Huxley-type equation with $K=1$ and $s=1$.}
\label{fig:5} 
\end{figure}

The exact solution for the population is shown in Figure \ref{fig:5}, which shows the population density decreasing over time. In this case, in general  the diffusivity is not a monotonic function of density but it has a local maximum. At late times, $\theta$ is small, $D(\theta)$ is increasing for small $\theta$, and the shape of the solution approaches the $J_0$ Bessel function from below. At early times, the range of densities may allow $D(\theta)$ to be non-monotonic, and the scaled density may lie above the limiting Bessel function.

In the population genetics application, this means that a new recessive engineered gene can be completely removed from an isolated population simply by extreme culling at the boundary. In the application to the nerve axon potential, it means that the potential can be set to a uniform value internally by applying that potential at the boundary.

\section{Fitzhugh-Nagumo equation, with $R(\theta)=s\theta(1-\theta)(\theta-\theta_1)$}

The cubic reaction term as an appropriate model in biological applications was first introduced in the late 1960's and early 1970s (see for example \cite{N75}). In this case, the generalisation to higher dimensions and to include nonlinear diffusion is
\[
\theta_t=\nabla\cdot[D(\theta)\nabla\theta]+s\theta(1-\theta)(\theta-\theta_1).
\]

Whereas the Fisher and Huxley equations have a single stable steady state at $\theta=1$, the Fitzhugh Nagumo equation is a bi-stable model, with an intermediate unstable steady state $\theta=\theta_1$ between the stable activated state $\theta=1$ and stable quiescent state $\theta=0$. The model is best known in one spatial dimension for propagation of the axon potential. It could be considered in higher dimensions for bi-stable activation of a specific neurological function associated with one region of the nervous system. The model arises also in diploid population genetics in the generic case that three possible genotypes have fitness coefficients that are neither Mendelian nor in arithmetic progression \cite{BBFA02,BB04}. It could also represent a population that is artificially introduced to a new environment where it must have a threshold density of $\theta_1$ in order to survive.

First, consider the case when $\kappa=0$. The solution to equation \eqref{use} is
\[
D(\theta)=c_1\theta^{-(1+\alpha/\theta_1)}(\theta-1)^{-1+\alpha/(\theta_1-1)} (\theta-\theta_1)^{-1-\alpha/(\theta_1(\theta_1-1))},
\]
where $\alpha=A/s$. When $A<0$, the expression for $D(\theta)$ is divergent as $\theta\to\theta_1$ and when $A>0$, it is divergent at both $\theta=0$ and $\theta=1$. When $A=0$, we see that $D(\theta)$ is divergent at all three zeros of $R(\theta)$. As such, the case when $\kappa=0$ is inadmissible when modelling a population.

When $\kappa=-K^2<0$, we can use either the fact that $R(1)=0$ or $R(\theta_1)=0$ to show that we must have $A>0$ in order to have $D(\theta)$ positive. We may then use the same arguments presented in section 3 to conclude that this case is also inadmissible.

Finally, we consider the case $\kappa=K^2>0$. From \eqref{use}, we see that 
\[
D(0)=\lim_{\theta\to0}\ds\frac{AD(0)\theta}{s\theta(-\theta_1)-K^2D(0)\theta}=\ds\frac{AD(0)}{-s\theta_1-K^2D(0)}~,
\]
except in the singular case of the the denominator being zero, $D(0)=-(s\theta_1)/K^2$. Also $D(\theta_1)=D(1)=-A/K^2$. Since $D>0$, we can deduce that $A<0$ and $u$ will approach zero exponentially in time from \eqref{Helm}. Also, in the non-singular case,
\[
A=-(s\theta_1+K^2D(0)).
\]

This may represent an axon potential or a population decaying towards extinction,  because of culling at the boundary even if it started above the activation threshold.  For a 2-dimensional problem with rotational invariance, the solution to equation \eqref{Helm} is the Bessel function, $J_0(Kr)$. Once again, we can deduce the minimum size of an aquatic reserve such that the undesirable extinguishing solution does not exist. If $r_1$ is the radius of the aquatic reserve, outside of which the population is practically zero, then $r_1=\lambda_1/K$, where $\lambda_1$ is the first zero of the Bessel function. For the population to be maintained, we require $K^2=-s\theta_1/D(0)$. For a population with a carrying capacity at $\theta=1$, it must be true that $s>0$ so that the local growth rate is negative at higher populations. Since $K^2>0$ and $s>0$, this means that $\theta_1<0$. In that case, the quiescent state (or extinction state) is now the intermediate unstable steady state. It then follows that
\[
r_1>\lambda_1\sqrt{\frac{D(0)}{s|\theta_1|}}
\]
is the minimum radius for a marine park in order for the extinction point not to be stable. For example and using the same parameters as before, if $D(0)=100$km$^2$ per year, $1/s=5$ years, and $|\theta_1|=0.4$, the reserve should be at least $2r_1=170$km in diameter.

Once again, an estimate for $D(\theta)$ can be calculated using the iterative map \eqref{Dmap}. Since the carrying capacity better represents the ecological potential, we use the value $D_0=D(1)=-A/K^2$ as a starting value for the iterative procedure, guaranteeing that all subsequent iterates have the correct value at $\theta=1$. The first two iterates are:
\begin{eqnarray*}
D_0&=&-\ds\frac{A}{K^2}\\
D_1&=&\ds\frac{A^2}{K^2s}\frac{1}{P(\theta)}
\end{eqnarray*}
where $P(\theta)=\theta^2-(\theta_1+1)\theta+\theta_1-A/s$. The form for $D_2$ will depend on the values of the parameters $\theta_1$, $A$ and $s$. We have the necessary condition that $\theta_1<|A|/s$ since we must have $D(0)>0$. If $P(\theta)$ has two real roots, then we also have that either $\theta_1>1+2\sqrt{|A|/s}$ or $\theta_1<1-2\sqrt{|A|/s}$. In this case, the next iterate, $D_2$ may be written in terms of $\log$ functions and, as such, it becomes singular for certain values of $\theta$ in the domain of interest.

If $P(\theta)$ has a repeated root, the next iterate will depend inversely upon a quartic in $\theta$. This iterate for the diffusivity then also has singularities for certain values of $\theta$in the domain of interest and as such, we proceed no further.

\begin{figure}
\includegraphics[trim = 15mm 70mm 15mm 70mm, clip, width=100mm]{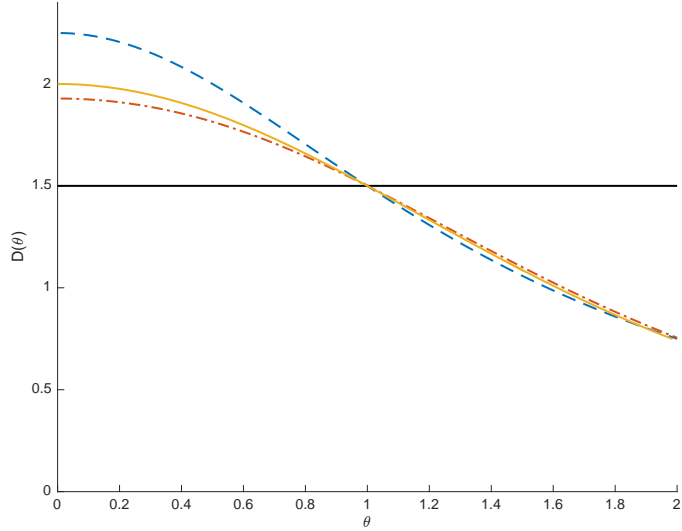}
\caption{$D(\theta)$ for the Fitzhugh-Nagumo source term, constructed by \mat routine ode45 (solid), as well as approximations $D_0$ (solid), $D_1(\theta)$ (dashed), $D_2(\theta)$ (dash-dot) ($K=1$, $s=1$ and $A=-1.5$).}
\label{fig:6} 
\end{figure}

If $P(\theta)$ has two complex roots, then we require $1-2\sqrt{|A|/s}<\theta<1+2\sqrt{|A|/s}$ as well as $\theta_1<|A|/s$. In this case, the integral of $D_1$ is an inverse tan function, and the next iterate, $D_2$, is given by
\[
D_2=\ds\frac{A^3}{K^2}~\ds\frac{\atan\left(\frac{1}{\beta}\left(\theta-\nfrac12(\theta_1+1)\right)\right)+\atan\left(\frac{1}{2\beta}(\theta_1+1)\right)}{\beta s^2\theta(1-\theta)(\theta-\theta_1)- A^2\left[\atan\left(\frac{1}{\beta}\left(\theta-\nfrac12(\theta_1+1)\right)\right)+\atan\left(\frac{1}{2\beta}(\theta_1+1)\right)\right]},
\]
where $\beta^2=-\nfrac14\theta_1^2+\nfrac12\theta_1-\nfrac14-A/s$. Some values of the parameters produce reasonable functions for $D_2$, for example, if $K=1$, $A=-1.5$, $s=0.5$ and $\theta_1=-1$, the first three iterates are shown in Figure \ref{fig:6}, together with the numerical solution to equation \eqref{use} found using \mat's ode45 routine.

\begin{figure}
\includegraphics[trim = 35mm 100mm 35mm 100mm, clip, width=100mm]{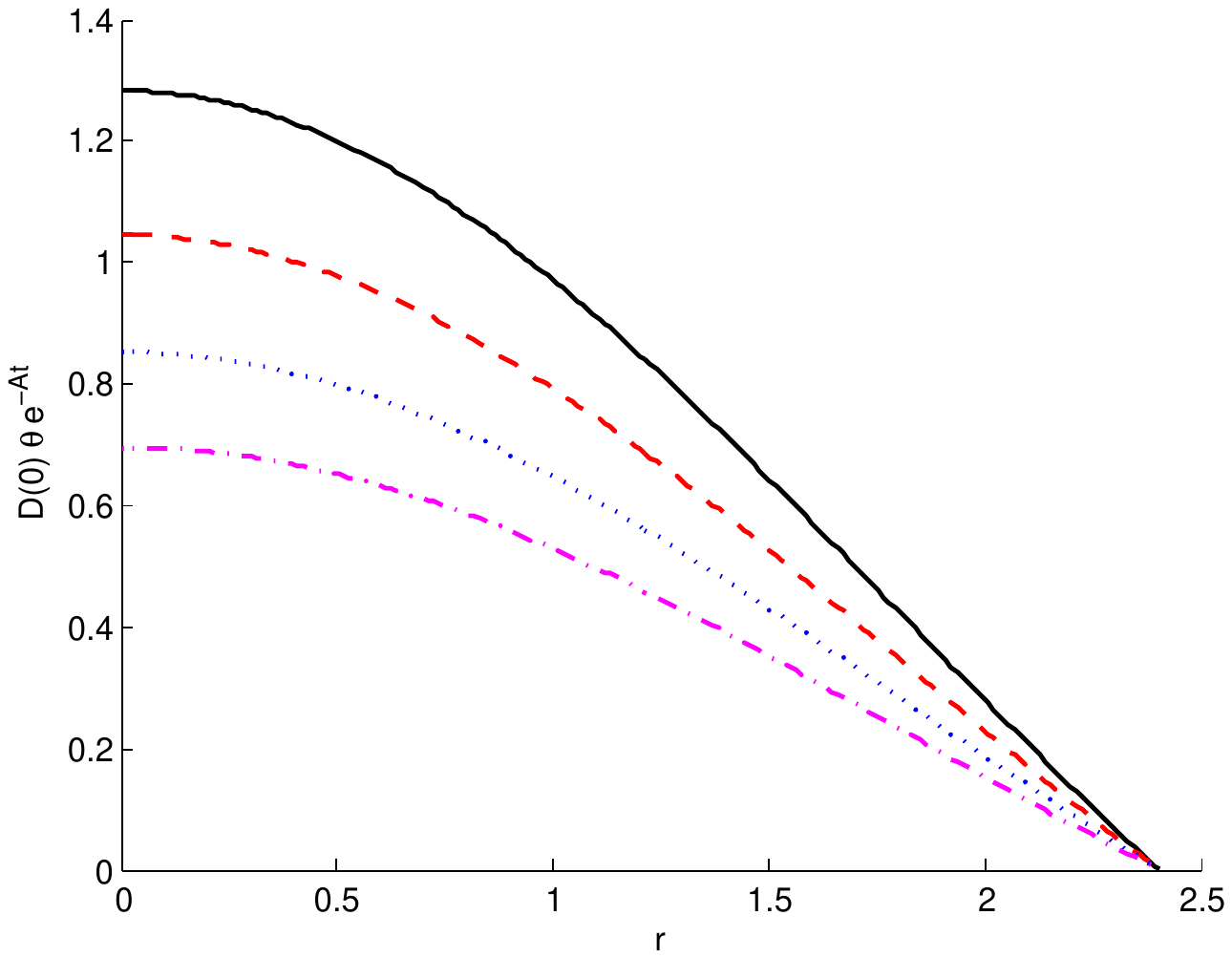}
\caption{Solution for the decreasing population density versus radial coordinate, for the Fitzhugh-Nagumo source term, at times $|A|t=0.0, 0.1, 0.2, 0.3$, according to the Fitzhugh-Nagumo type equation with $K=1$ and $s=1$.}
\label{fig:7} 
\end{figure}

The exact solution for the population is plotted in Figure \ref{fig:7}, showing the population is doomed to extinction.

\section{Application to population genetics}

Let $\theta$ be the frequency of new advantageous allele $a_2$, and $1-\theta$ be the frequency of the original allele $a_1$. Let $\gamma_{11}, \gamma_{12}, \gamma_{22}$ be the fitness coefficients of genotypes with zero, one and two copies of $a_2$ respectively. Allele $a_2$ will be said to be partially recessive if $\gamma_{22}-\gamma_{12}>\gamma_{12}-\gamma_{11}$. That is, addition of the second copy of gene $a_1$ gives a greater advantage than does the first copy. For example, in the reaction diffusion equations with linear diffusion, calculations in \cite{BB04} show that the generic source term is that of Fitzhugh-Nagumo-type, $s\theta (1-\theta)(\theta-\theta_1)$ with 
$$
s=\gamma_{11}-2\gamma_{12}+\gamma_{22};\quad\quad\theta_1=\frac{1}{2-\nu};\quad\quad\nu=\frac{\gamma_{22}-\gamma_{11}}{\gamma_{22}-\gamma_{12}}
$$

For a partially recessive gene $a_2$, $\nu>2$ and $\theta_1 <0$. For the nonlinear diffusion model with Fitzhugh-Nagumo-type source, the extinguishing solution developed in Section 5 exists only if the radius of the boundary where the new gene is selectively removed, is less than the critical value $r_c=\lambda_1\sqrt{D(0)/s|\theta_1|}$.

The case $\nu=1$ is a special case that gives rise to the Fisher equation with quadratic logistic source term $s\theta(1-\theta)$. This arises in the context of diploid sexual population genetics only in the special case that the fitness coefficients are in arithmetic progression, ie $\gamma_{22}-\gamma_{12}=\gamma_{12}-\gamma_{11}$. In that case, allele $a_2$ is neither partially recessive nor partially dominant.  The Fisher equation arises also as the equation for the frequency of a new advantageous gene of an asexual species \cite{BBFA02}. The critical radius for existence of the extinguishing solution in the case of the quadratic logistic source term is $r_c=\lambda_1\sqrt{D(0)/s}$.

Another special case arises when $\gamma_{11}=\gamma_{12}<\gamma_{22}$, the case of a totally recessive gene $a_2$. The phenotype whose expression requires two copies of a new gene, is particularly vulnerable.  Now $\theta_1=1$, the value at which the critical radius $r_c$ diverges. The extinguishing solution of the Huxley reaction-diffusion equation (Section 4) always exists, no matter how large is the diameter of the boundary where individuals with the totally recessive new gene are removed.

If a new, perhaps genetically engineered, advantageous gene is totally recessive or if it is partially recessive within a population that is contained within a circle smaller than that with critical radius, then the frequency of the new gene can be reduced to zero uniformly just by selective culling at the boundary.

If the new gene is partially recessive and the circular domain of the species has a radius larger than the critical value, then its removal cannot be achieved by actions taken at the boundary alone.

\section{Discussion and final remarks}

Reaction-diffusion equations, with a source term that is either quadratic or cubic, are commonly used to model various physical and biological systems. Here we have shown that equations of this type with nonlinear diffusion, are exactly solvable provided the nonlinear diffusion and source terms satisfy a particular relationship. This relationship can be found by first writing the equation in terms of the Kirchhoff variable, and then using the nonclassical symmetry reduction method. Provided that the relationship is satisfied, the equation may be separated into spatial and temporal parts. The solutions for the Kirchhoff variable are then exponential in time (either growth or decay), and satisfy either a Helmholtz, modified Helmholtz or Laplace equation in space.

In general, one may either specify the nonlinear diffusivity and then calculate the corresponding reaction term which will allow exact solution, or alternatively, the reaction term may be specified and the diffusivity calculated using the aforementioned relationship. In this paper, we are particularly interested in logistic Fisher-type models and cubic Huxley- or Fitzhugh Nagumo-type models, and so we calculate the nonlinear diffusivity after specifying the reaction term. In this case, the diffusivity can be constructed by solving a differential equation that is equivalent to an Abel equation if the Kirchhoff variable satisfies a Helmholtz (or modified Helmholtz) equation, or a separable equation if the Kirchhoff variable satisfies the Laplace equation. In cases where the nonlinear diffusivity cannot be written in closed form, it may be written in terms of a series expansion. It is important to note that, although the diffusivity is required to be a nonlinear function of the dependent variable, for reasonable parameters sets, it is well behaved and slowly varying over the relevant range of the population density.

Interpretation of the obtained solutions, in the application of population genetics, is discussed in Section 6.  The extinguishing solution of the Huxley reaction-diffusion equation (Section 4) always exists, no matter how large is the diameter of the boundary where individuals with the totally recessive new gene are removed. In practice, it is more likely that only individuals with the distinctive phenotype expressed from two copies of the recessive gene would be identified and removed at the boundary. This would locally remove a fraction $[2\theta^2/(2\theta^2+2\theta[1-\theta])]=\theta$ of the new genes from the gene pool. This suggests a nonlinear Robin boundary condition $-u_r=h(u)\theta$, at $r=r_1$, representing a selective removal from the domain of the new gene, in proportion to its frequency. If only the proportion $\theta$ of the new genes are detected, then the radiation constant $h(u(\theta))$ is proportional to $\theta$, $h=h_1\theta$, which is approximately $h_1u/D(0)$ with $h_1$ constant. The nonlinear radiation boundary condition therefore is 
$$
-u_r=Hu^2 ;~~r=r_1,
$$
with $H=h_1/D(0)$ (constant). This boundary condition would be of interest for the application to population genetics.

\end{document}